\documentclass{amsart}
\usepackage{amsmath,amsthm,amssymb,comment}
\usepackage[pagewise]{lineno}
\theoremstyle{definition}
\newtheorem{definition}{Definition}[section]
\theoremstyle{plain}

\newtheorem{theorem}[definition]{Theorem}
\newtheorem{corollary}[definition]{Corollary}
\theoremstyle{remark}
\newtheorem{remark}[definition]{Remark}

\DeclareMathOperator{\End}{End}

\title{Weierstrass-type representations}
\author{Mason Pember}

\begin{document}
\maketitle
\begin{abstract}
Weierstrass-type representations have been used extensively in surface theory to create surfaces with special curvature properties. In this paper we give a unified description of these representations in terms of classical transformation theory of $\Omega$-surfaces. 
\end{abstract}

\section{Introduction}
The Weierstrass-Enneper representation~\cite{W1866} is a fundamental tool in differential geometry for creating interesting examples of minimal surfaces in Euclidean 3-space. A similar representation was developed by Bryant~\cite{B1987}, and later~\cite{UY1992}, for surfaces of constant mean curvature (CMC) $H\equiv c$ in hyperbolic 3-space of constant sectional curvature $-c^{2}$. Many other Weierstrass-type representations exist for various surface classes, for example: 
\begin{itemize}
\item maximal surfaces in Lorentz 3-space~\cite{K1983},
\item CMC surfaces $H\equiv c$ in de Sitter 3-space of constant sectional curvature $c^{2}$~\cite{AA1998},
\item flat surfaces in hyperbolic 3-space~\cite{GMM2000}, 
\item linear Weingarten surfaces of Bryant type in hyperbolic 3-space~\cite{GMM2004}, 
\item linear Weingarten surfaces of Bianchi type in de Sitter 3-space~\cite{AE2007}.
\end{itemize}
The ingredients for such representations are always the same: a meromorphic function and a holomorphic 1-form. Thus one might expect that these representations are related in some way. A unification of some of these representations was achieved in~\cite{AGM2005} under the umbrella of marginally trapped surfaces in Minkowski space. 

On the other hand, geometric interpretations of these representations have been sought in various works. The classical Weierstrass-Enneper representation can be understood using the Christoffel transformation of isothermic surfaces~\cite{H2003}. The Umehara-Yamada perturbation~\cite{UY1992} deforms minimal surfaces in Euclidean 3-space into CMC-1 surfaces in hyperbolic 3-space. This perturbation was given a M\"{o}bius geometric interpretation in~\cite{HMN2001}. A Laguerre geometric interpretation of this perturbation was studied in~\cite{MN2016}, extending this notion to a wider class of surfaces.  
 
Motivated by an observation in~\cite{BHR2012}, we seek to understand Weierstrass-type representations using transformations of $\Omega$-surfaces. In~\cite{D1911iii}, using the Christoffel transformation for isothermic surfaces in Minkowski space, Demoulin develops a notion of dual surfaces for $\Omega$-surfaces, yielding a transformation for this surface class. We show that L-isothermic surfaces can be characterised as 
$\Omega$-surfaces admitting a special dual surface, akin to how minimal surfaces in Euclidean space can be characterised as Christoffel transformations of their Gauss map. We then give an invariant explanation of all of the aforementioned Weierstrass-type representations as an application of Demoulin's dual transformation to a prescribed Gauss map.  

\textit{Acknowledgements.} The author would like to thank F. Burstall, J. Cho, U. Hertrich-Jeromin, W. Rossmann and M. Yasumoto for fruitful and pleasant conversations about this topic. This work was supported by the Austrian Science Fund (FWF) through the research project P28427-N35 ``Non-rigidity and symmetry breaking".

\section{$\Omega$-surfaces}

Let $x:\Sigma\to \mathbb{R}^{3,1}$ be a spacelike immersion. As defined in~\cite{C1867,P1988} a smooth map $x^{*}:\Sigma\to \mathbb{R}^{3,1}$ is called a \textit{Christoffel dual} of $x$ if $x$ and $x^{*}$
\begin{itemize}
\item have parallel tangent planes, 
\item induce the same conformal structure on $T\Sigma$, and
\item induce opposite orientations on $T\Sigma$.
\end{itemize}
It was shown\footnote{In fact, it was shown using Clifford algebra that this relation can be expressed by the vanishing of a single wedge.} in~\cite{B2006} that these conditions can be encapsulated by
\begin{equation}
\label{eqn:cdual}
(dx\wedge dx^{*})=0 \quad \text{and}\quad dx\curlywedge dx^{*}=0, 
\end{equation}
where $(dx\wedge dx^{*})$ is a symmetric 2-form defined by 
\[ (dx\wedge dx^{*})(X,Y) = (dx(X),dx^{*}(Y))-(dx(Y),dx^{*}(X)),\]
and $dx\curlywedge dx^{*}$ is a symmetric 2-form with values in $\wedge^{2}\underline{\mathbb{R}}^{3,1}$ (with $\underline{\mathbb{R}}^{3,1}$ denoting the trivial bundle $\Sigma\times\mathbb{R}^{3,1}$) defined by 
\[ dx\curlywedge dx^{*}(X,Y) = dx(X)\wedge dx^{*}(Y)-dx(Y)\wedge dx^{*}(X).\]
A surface which possesses a Christoffel dual is then said to be \textit{isothermic}. By the symmetric nature of~\eqref{eqn:cdual}, $x^{*}$ is itself isothermic with Christoffel dual $x$. 

In Laguerre geometry, one uses isotropy projection to identify points in $\mathbb{R}^{3,1}$ with spheres (see for example~\cite{C2008}). Therefore a spacelike immersion $x:\Sigma\to\mathbb{R}^{3,1}$ represents a congruence of spheres. Given a spacelike immersion $x:\Sigma\to\mathbb{R}^{3,1}$, we may write the normal bundle of $x$ as $dx(T\Sigma)^{\perp}=G\oplus \tilde{G}$, where $G$ and $\tilde{G}$ are null line subbundles of $dx(T\Sigma)^{\perp}$. We call $G$ and $\tilde{G}$ the \textit{lightlike Gauss maps} of $x$. We may then construct null affine line bundles $L:=x+G$ and $\tilde{L}:=x+\tilde{G}$ of $\underline{\mathbb{R}}^{3,1}$. In Laguerre geometry, these represent the envelopes of the sphere congruence $x$ (see for example~\cite{MN2000}). 

As classically defined by Demoulin~\cite{D1911iii,D1911i,D1911ii}, the envelopes of isothermic (spacelike) sphere congruences are called $\textit{$\Omega$-surfaces}$. Therefore, $\Omega$-surfaces are the null affine line bundles $L$ of $\underline{\mathbb{R}}^{3,1}$, for which we may write $L=x+G$ for some isothermic surface $x:\Sigma\to\mathbb{R}^{3,1}$. Given Christoffel dual maps $x$ and $x^{*}$, we say that envelopes $L$ and $L^{*}$ of $x$ and $x^{*}$, respectively, are \textit{$\Omega$-dual} if $L$ is parallel\footnote{It is always possible to arrange $L$ and $L^{*}$ to be parallel, since the normal bundles of $x$ and $x^{*}$ are the same.} to $L^{*}$. 

\section{Marginally trapped surfaces}
\label{sec:margtrap}

A spacelike immersion $x:\Sigma\to \mathbb{R}^{3,1}$ is called \textit{marginally trapped} if its mean curvature vector
\[ \textbf{H} = \frac{1}{2}(d_{X}d_{X}x + d_{Y}d_{Y}x)^{\perp}\]
is lightlike, where $X,Y\in \Gamma T\Sigma$ is an orthonormal basis with respect to the induced metric of $x$ and $(.)^{\perp}$ denotes projection onto the normal bundle of $x$. Given $p\in \Sigma$, $\textbf{H}(p)$ is then lightlike if and only if $\textbf{H}(p)\in G(p)$ or $\textbf{H}(p)\in \tilde{G}(p)$, where $G$ and $\tilde{G}$ are the lightlike Gauss maps of $x$. We say that a marginally trapped surface $x$ is \textit{regular} if either $\textbf{H}(p)\in G(p)$ for all $p\in\Sigma$ or $\textbf{H}(p)\in \tilde{G}(p)$ for all $p\in\Sigma$, i.e., $\textbf{H}\in \Gamma G$ or $\textbf{H}\in \Gamma \tilde{G}$. 

Assume now that the normal bundle of $x$ is flat. Then there locally exists parallel sections $g\in \Gamma G$ and $\tilde{g}\in \Gamma \tilde{G}$ with $(g,\tilde{g})=-1$. Since $g$ is parallel, we may write $dg=dx\circ S$ for some $S\in \Gamma \End(T\Sigma)$. Since $(dx,dg)$ is symmetric, we have that $S$ is symmetric with respect to the induced metric of $x$. Therefore, $S$ admits a basis $X,Y\in \Gamma T\Sigma$ of orthonormal (with respect to the induced metric of $x$) eigenvectors with respective eigenvalues $\alpha$ and $\beta$. Now the mean curvature vector field is given by 
\begin{align*}  
\textbf{H}=- \frac{1}{2}((d_{X}d_{X}x + d_{Y}d_{Y}x,g)\tilde{g} + (d_{X}d_{X}x + d_{Y}d_{Y}x,\tilde{g})g)= \frac{1}{2}(\alpha + \beta)\tilde{g} \, mod\, G.
\end{align*}
Thus, $\textbf{H}\in \Gamma G$ if and only if $\alpha+\beta=0$, i.e., $S$ is trace-free. On the other hand 
\[ dx\curlywedge dg(X,Y) = (\beta+ \alpha) d_{X}x\wedge d_{Y}x.\]
Thus $S$ is trace-free if and only if $dx\curlywedge dg=0$. 

If $dx\curlywedge dg=0$ for some $g\in \Gamma G$, then $(ddx,g) =- (dx\wedge dg)=0$. Thus, $dx\curlywedge dg=0$ implies that $x$ and $g$ are Christoffel dual. We thus arrive at the following theorem:

\begin{theorem}
A spacelike immersion in Minkowski space is a regular marginally trapped surface with flat normal bundle if and only if it is Christoffel dual to a section of one of its lightlike Gauss maps. 
\end{theorem}

Suppose that $x:\Sigma\to\mathbb{R}^{3,1}$ is Christoffel dual to $g\in \Gamma G$. By defining $L:=x+G$, we see that $L$ is $\Omega$-dual to $L^{*}:=G$. Moreover, since $x$ has the same induced conformal structure on $T\Sigma$ as $G$, one identifies $x$ as the \textit{middle sphere congruence} of $L$, see~\cite{B1929,MN2018}. One deduces from~\cite{BHPR2018} that such $L$ are the envelopes of \textit{L-isothermic surfaces}, that is, surfaces in Euclidean 3-space that admit curvature line coordinates which are conformal with respect to the third fundamental form. This characterisation gives an analogue of minimal surfaces being Christoffel dual to their Gauss map in Euclidean space: 

\begin{corollary}
Envelopes of L-isothermic surfaces are the $\Omega$-surfaces that are $\Omega$-dual to their lightlike Gauss map. 
\end{corollary}

Now, the condition $dx\curlywedge dg =0$, implies that\footnote{Throughout this paper we shall use the well-known identification $\wedge^{2}\mathbb{R}^{3,1}\cong \mathfrak{so}(3,1)$, via $(a\wedge b)v = (a,v)b-(b,v)a$.} $\zeta:= g\wedge dx\in \Omega^{1}(G\wedge G^{\perp})$ is a closed 1-form. Moreover, since $G\wedge G^{\perp}$ is an abelian subbundle of $\wedge^{2}\underline{\mathbb{R}}^{3,1}$, one has that $[\zeta\wedge \zeta]=0$. Therefore, $\{d+t\zeta\}_{t\in\mathbb {R}}$ is a 1-parameter family of flat metric connections. Let $T_{t}:\Sigma\to O(3,1)$ be the local orthogonal trivialising gauge transformations of $d+t\zeta$, i.e., 
\[ T_{t}(d+t\zeta)T_{t}^{-1} = d.\]
Such transformations are unique up to premultiplication by a constant endomorphism $A\in O(3,1)$. Now $T_{t}dx\in \Omega^{1}(\underline{\mathbb{R}}^{3,1})$ is a closed 1-form: 
\[ d(T_{t}dx) = T_{t}(d+t\zeta)dx = T_{t}(ddx + t\zeta\wedge dx) =0\]
since $\zeta\wedge dx =-(dx\wedge dx) g = 0$. Thus we may integrate to obtain a new surface $x_{t}:\Sigma\to\mathbb{R}^{3,1}$ satisfying $dx^{t} = T_{t}dx$. Define $G_{t}:= T_{t}G$ and set $g_{t}:=T_{t}g$. Then $dg_{t} = T_{t}(d+t\zeta)g=T_{t}dg$, since $\zeta g=0$, and one deduces that 
\[ dx_{t}\curlywedge dg_{t} = T_{t} \, dx\curlywedge dg \, T_{t}^{-1} = 0.\]
Thus $x^{t}$ is a marginally trapped surface and $L_{t}:= x_{t}+G_{t}$ is the envelope of an L-isothermic surface. This is the \textit{T-transform} of L-isothermic surfaces (see~\cite{MN2016,MN2018}). Notice that we obtain a new closed 1-form 
\[ \zeta_{t}:= g_{t}\wedge dx_{t} = Ad_{T_{t}}\zeta\in \Omega^{1}(G_{t}\wedge G_{t}^{\perp}) .\]
We therefore obtain local orthogonalising gauge transformations $T^{t}_{s}:\Sigma\to O(3,1)$ of the 1-parameter family of flat connections $\{d+s\zeta_{t}\}_{s\in\mathbb{R}}$. In~\cite[Section 5.5.9]{H2003} the following property was shown for iterating these transformations: 
\begin{equation}
\label{eqn:iterate}
T^{t}_{s}T_{t} = T_{t+s},\quad \text{and thus}\quad  T^{t}_{-t} = (T_{t})^{-1}. 
\end{equation}
This property will be useful for us in the following section. 

\section{Weierstrass-type representations}

The ingredients for Weierstrass-type representations are:
\begin{itemize}
\item \textbf{A simply connected Riemann surface $\Sigma$}. Equivalently, $\Sigma$ is a simply connected 2-dimensional manifold 
equipped with a conformal structure and an orientation. 

\item \textbf{A meromorphic function $\phi:\Sigma\to \mathbb{C}\cup \{\infty\}$}. Equivalently, since $\mathbb{C}\cup \{\infty\} \cong S^{2}\cong \mathbb{P}(\mathcal{L})$, where $\mathcal{L}\subset\mathbb{R}^{3,1}$ denotes the lightcone, we can identify $\phi$ with a smooth map $G:\Sigma\to \mathbb{P}(\mathcal{L})$. $\phi$ being meromorphic is equivalent to $G$ being an orientation preserving map whose induced conformal structure is weakly equivalent to the conformal structure on $\Sigma$. 

\item \textbf{A holomorphic 1-form $\omega$}. Alternatively, one may prescribe a holomorphic quadratic differential\footnote{That is, we may write $q^{2,0} = h dz^{2}$ for some local holomorphic coordinate $z$ on $\Sigma$ and some holomorphic function $h$.} $q$. We then have the relation $q=d\phi\, \omega + d\bar{\phi}\,\bar{\omega}$. 
\end{itemize}
We make the assumption\footnote{Since $\phi$ is meromorphic, this excludes the case that $G$ is constant and isolated points of $\Sigma$ where $G$ does not immerse.} that $\phi$ has no critical points, and thus $G$ is an immersion. Now for any non-zero lift $g\in \Gamma G$, we may write $q = (dg,dg\circ Q)$ for some endomorphism $Q\in \Gamma \End(T\Sigma)$. $q$ is then a holomorphic quadratic differential if and only if
\[ \zeta := g\wedge dg\circ Q\in \Omega^{1}(G\wedge G^{\perp})\]
is a closed 1-form (see~\cite{BS2012,S2008}). 

Explicitly, one may identify a meromorphic function $\phi$ with the map $G:\Sigma\to\mathbb{P}(\mathcal{L})$ spanned by the lift\begin{equation}
\label{eqn:secG}
 g=(1+\phi\bar{\phi})e_{0} + (\phi +\bar{\phi})e_{1} - i(\phi - \bar{\phi}) e_{2} + (-1+\phi\bar{\phi})e_{3}\in \Gamma G,
\end{equation}
where $\{e_{0},...,e_{3}\}$ is a pseudo-orthonormal basis for $\mathbb{R}^{3,1}$ with $e_{0}$ timelike and $e_{1},e_{2},e_{3}$ spacelike. The induced metric of $g$ is $(dg,dg) = 4 d\phi d\bar{\phi}$. By defining 
\[ Q = \frac{1}{2}\left( \bar{\omega}\otimes \frac{\partial}{\partial \phi}  +  \omega\otimes\frac{\partial}{\partial \bar{\phi}}\right),\]
we have that $(dg,dg\circ Q) = d\phi \omega + d\bar{\phi}\bar{\omega} =q$. One then computes
\begin{align}
\label{eqn:zeta}
\zeta = g\wedge dg\circ Q &= \frac{1}{2}\{(e_{0}+\phi e_{1}-i\phi e_{2}-e_{3})\wedge (\phi e_{0} + e_{1} + i e_{2} + \phi e_{3})\omega\\
&+ (e_{0}+\bar{\phi} e_{1}+i\bar{\phi} e_{2}-e_{3})\wedge (\bar{\phi} e_{0} + e_{1} - i e_{2} + \bar{\phi} e_{3})\bar{\omega}\}\nonumber.
\end{align}
Notice that given Weierstrass data $(\phi, \omega)$, the right hand side of~\eqref{eqn:zeta} yields a closed 1-form with values in $G\wedge G^{\perp}$, regardless of whether $\phi$ has critical points or not. 

\subsection{Affine hyperplanes in $\mathbb{R}^{3,1}$}
\label{subsec:affine}
Let $\mathfrak{p}\in\mathbb{R}^{3,1}$ be a non-zero vector. Then hyperplanes with normal $\mathfrak{p}$ are flat 3-dimensional affine spaces. Given a closed $\zeta\in \Omega^{1}(G\wedge G^{\perp})$ we have that $\zeta\mathfrak{p}\in \Omega^{1}(\underline{\mathbb{R}}^{3,1})$ is a closed 1-form and we may locally integrate it to obtain 
\[ x:\Sigma\to \mathbb{R}^{3,1}\quad \text{satisfying} \quad dx = - \zeta\mathfrak{p}.\] 
Now $d(x,\mathfrak{p})= (\zeta\mathfrak{p},\mathfrak{p})=0$, since $\zeta$ is skew-symmetric. Thus $x$ takes values in an affine hyperplane with normal $\mathfrak{p}$. If $G(p)\perp \mathfrak{p}$ for some $p\in \Sigma$ then $d_{p}x\in G(p)$ and thus $x$ does not immerse at $p$. Away from such points, we may write $\zeta = g\wedge \omega$, where $g\in \Gamma G$ satisfies $(g,\mathfrak{p})=-1$ and $\omega\in \Omega^{1}(G^{\perp})$. Then $dx=-\zeta\mathfrak{p} = \omega \, mod\, G$ implies that $\zeta = g\wedge dx$ and the closedness of $\zeta$ implies that $dg\curlywedge dx=0$. Hence, $L:=x+G$ is $\Omega$-dual to $L^{*}:= G$ and $x$ is marginally trapped. In fact, since $\mathfrak{p}$ lies in $dx(T\Sigma)^{\perp}$, one deduces that $\textbf{H}=0$. Hence, $x$ has zero mean curvature. We thus have the following 3 cases:
\begin{enumerate}
\item if $\mathfrak{p}$ is timelike then $x$ is a minimal surface in a Euclidean 3-space,
\item if $\mathfrak{p}$ is spacelike then $x$ is a maximal surface in a Lorentzian 3-space, 
\item if $\mathfrak{p}$ is lightlike then $x$ has zero mean curvature in an isotropic 3-space. 
\end{enumerate}
In cases (1) and (2) we obtain a unit (spacelike or timelike, respectively) normal of $x$ by setting $n:= g-\mathfrak{p} \in \Gamma \mathfrak{p}^{\perp}$. 

By choosing $\mathfrak{p}=e_{0}$, one deduces from~\eqref{eqn:zeta} that 
\[ dx = -\zeta e_{0} = \text{Re}\{((1- \phi^{2})e_{1} + i(1 + \phi^{2})e_{2} +2\phi e_{3})\omega\}\]
and we thus recover the Weierstrass-Enneper representation~\cite{W1866}. By choosing $\mathfrak{p}=e_{3}$, we have that 
\[ dx = - \zeta e_{3} = \text{Re}\{(2\phi e_{0} + (1+ \phi^{2})e_{1} + i(1 -\phi^{2})e_{2})\omega\},\]
recovering the representation of~\cite{K1983} for maximal surfaces in Minkowski 3-space. Choosing $\mathfrak{p}= \frac{e_{0}+e_{3}}{2}$ we obtain a representation of zero mean curvature surfaces in isotropic 3-space:
\[ dx = -\zeta\tfrac{e_{0}+e_{3}}{2}= \text{Re}\left\{\left( e_{1}+ie_{2} + 2\phi\tfrac{e_{0}+e_{3}}{2}\right)\omega\right\}.\]
\begin{remark}
In~\cite{MN2016} it was shown that surfaces in Euclidean space that are simultaneously L-minimal and L-isothermic are those surfaces whose middle sphere congruence is one of the three cases above. On the other hand, a Weierstrass-type representation was developed for such surfaces in~\cite{S2009}. One can recover this representation by intersecting the envelopes $L$ of the three cases above with appropriate affine Euclidean 3-spaces. 
\end{remark}

\subsection{Quadrics in $\mathbb{R}^{3,1}$}
\label{subsec:quadrics}
Given a closed 1-form $\zeta\in \Omega^{1}(G\wedge G^{\perp})$ we have that $\{d+t\zeta\}_{t\in\mathbb{R}}$ is a 1-parameter family of flat connections. Therefore, there locally exists parallel sections of these connections. Suppose that $x:\Sigma\to \mathbb{R}^{3,1}$ satisfies 
\[ (d+m\zeta)x = 0\]
for some\footnote{The choice of $m$ here amounts to a constant scaling of the Hopf differential. In many works this scaling is fixed by choosing $m$ appropriately.} non-zero $m\in\mathbb{R}$. Then $d(x,x) = -2m(\zeta x,x)=0$, since $\zeta$ is skew-symmetric. Hence, $(x,x)$ is constant. If $x(p)\perp G(p)$ for some $p\in \Sigma$, then $d_{p}x\in G(p)$ and thus $x$ does not immerse at $p$. Away from these points, we may write $\zeta = g\wedge \omega$ where $g\in \Gamma G$ such that $(g,x)=-1$ and $\omega\in \Omega^{1}(G^{\perp})$. Then the condition $dx= -m\zeta x$ implies that $\zeta = \frac{1}{m}g\wedge dx$. The closedness of $\zeta$ then implies that $dx\curlywedge dg=0$. Hence, $L:=x+G$ is $\Omega$-dual to $L^{*}=G$ and $x$ is marginally trapped. From~\cite{HI2015} we then obtain the following 3 cases:
\begin{enumerate}
\item if $(x,x)=-c^{2}$ then $x$ is a CMC-c surface in $\mathbb{H}^{3}(-c^{2})$,
\item if $(x,x)=c^{2}$ then $x$ is a CMC-c surface in $\mathbb{S}^{2,1}(c^{2})$,
\item if $(x,x)=0$ then $x$ is an intrinsically flat surface in $\mathcal{L}$. 
\end{enumerate}
In cases (1) and (2) $G$ then has a geometric interpretation as the hyperbolic Gauss map of $x$. An important observation is that parallel sections $x$ of $d+m\zeta$ are given by $x = T^{-1}_{m}\mathfrak{c}$, where $\mathfrak{c}\in \mathbb{R}^{3,1}$ and $T_{m}$ is a local trivialising orthogonal gauge transformations of $d+m\zeta$, i.e., $T_{m}(d+m\zeta)T_{m}^{-1}=d$ (see Section~\ref{sec:margtrap}). 

The Hermitian model of $\mathbb{R}^{3,1}$ identifies points in $\mathbb{R}^{3,1}$ with $2\times 2$ Hermitian matrices via the isometry
\[ x_{0}e_{0} + x_{1}e_{1} + x_{2} e_{2} + x_{3}e_{3} \mapsto \begin{pmatrix} x_{0}+ x_{3} & x_{1} + ix_{2}\\ x_{1}-ix_{2} & x_{0}-x_{3}\end{pmatrix},\]
where the metric on the space of Hermitian matrices is given by $(A,A)=-\det A$. Our basis vectors become
\[ e_{0} = \begin{pmatrix} 1&0\\0&1\end{pmatrix}, \, e_{1} = \begin{pmatrix} 0&1\\1&0\end{pmatrix},\, e_{2} = \begin{pmatrix} 0&i\\-i&0\end{pmatrix},\, e_{3}=\begin{pmatrix} 1&0\\0&-1\end{pmatrix}.\]
We identify $SL(2,\mathbb{C})$ with the orthogonal group $O(3,1)$ by its action on Hermitian matrices
\[ A\cdot v = AvA^{*}.\]
Skew-symmetric endomorphisms, i.e., elements of $\mathfrak{o}(3,1)$, are then identified with elements of $\mathfrak{sl}(2,\mathbb{C})$ via 
\[ B\cdot v = Bv + vB^{*}.\]
In particular, the skew-symmetric endomorphisms $e_{i}\wedge e_{j}$ are identified with $e_{ij}\in \mathfrak{sl}(2,\mathbb{C})$, where 
\[ e_{01} = -\frac{1}{2}e_{1}, \, e_{02} = -\frac{1}{2}e_{2}, \, e_{03} = -\frac{1}{2}e_{3}, \, e_{12} = \frac{i}{2} e_{3}, \, e_{13} = -\frac{i}{2}e_{2}, \, e_{23} = \frac{i}{2}e_{1}.\]
One then computes that $\zeta$ from~\eqref{eqn:zeta} is identified with 
\[ \xi = \frac{1}{2}\{ -(1-\phi^{2}) e_{1} - i(1+\phi^{2})e_{2} - 2\phi e_{3} \} \omega = \begin{pmatrix} -\phi & \phi^{2} \\ -1 &\phi\end{pmatrix}\omega \in \Omega^{1}(\mathfrak{sl}(2,\mathbb{C})).\]
Since $\xi$ is a closed 1-form, there exists $F_{t}:\Sigma\to SL(2,\mathbb{C})$ so that 
\[ F_{t}^{-1}dF_{t} = t\xi.\]
One quickly deduces that the orthogonal transformations identified with $F_{t}$ are in fact local trivialising orthogonal gauge transformations $T_{t}$ of $d+t\zeta$. Thus parallel sections of $d+m\zeta$ are given by $F^{-1}_{m}\mathfrak{c}(F^{-1}_{m})^{*}$ for $\mathfrak{c}\in\mathbb{R}^{3,1}$. Defining $\Psi_{m} := F^{-1}_{m}$ we have that 
\begin{equation}
\label{eqn:Psi}
d\Psi_{m} = -m\xi \Psi_{m} =m \begin{pmatrix} \phi & -\phi^{2} \\ 1 &-\phi\end{pmatrix}\omega \Psi_{m}.
\end{equation}
Now, by setting $\mathfrak{c} :=\frac{1}{2}((1-\mu)e_{0}+(1+\mu)e_{3})= \begin{pmatrix} 1 & 0\\ 0&-\mu\end{pmatrix}$ for $\mu\in \mathbb{R}$ and 
\begin{equation} 
\label{eqn:hyprep}
x = \Psi_{m}\mathfrak{c}(\Psi_{m})^{*}
\end{equation}
we obtain 
\begin{itemize}
\item CMC $H\equiv \frac{1}{\sqrt{-\mu}}$ surfaces in $\mathbb{H}^{3}(\frac{1}{\mu})$ when $\mu<0$; 
\item CMC $H\equiv \frac{1}{\sqrt{\mu}}$ surfaces in  $\mathbb{S}^{2,1}(\frac{1}{\mu})$ when $\mu>0$;
\item intrinsically flat surfaces in $\mathcal{L}$ when $\mu=0$. 
\end{itemize}
When $\mu = -1$ we see that this coincides with the representation of CMC $H\equiv 1$ in hyperbolic 3-space surfaces given in~\cite[Corollary 2.4]{RUY1997} and when $\mu=1$ this coincides with the representation of CMC $H\equiv 1$ in de-Sitter 3-space given in~\cite[Theorem 1]{F2007}. 

\subsection{The Umehara-Yamada perturbation}
Suppose that $x$ is a parallel section of $d+m\zeta$ for some non-zero $m\in \mathbb{R}$, and thus $x= T^{-1}_{m}\mathfrak{c}$ for some non-zero $\mathfrak{c}\in \mathbb{R}^{3,1}$. Then the T-transform $x_{m}$ of $x$ satisfies 
\[ dx_{m} = T_{m}dx= -T_{m}m\zeta x = - m\zeta_{m} T_{m}x = - m\zeta_{m}\mathfrak{c} .\] 
Hence, from Subsection~\ref{subsec:affine} we know that $x_{m}$ is a zero mean curvature surface in an affine hyperplane of $\mathbb{R}^{3,1}$. Therefore, we recover the result of~\cite{MN2016} that the T-transform of L-isothermic surfaces perturbs the cases (1), (2), (3) of Subsection~\ref{subsec:quadrics} into the cases (1), (2), (3) of Subsection~\ref{subsec:affine}, respectively. This generalises the Umehara-Yamada perturbation~\cite{UY1992} and gives a Laguerre geometric analogue of its interpretation in~\cite{HMN2001}. 

We call $G_{m}=T_{m}G$ the \textit{secondary Gauss map} of $x$. Since $G_{m}$ and $G$ induce the same conformal structures on $T\Sigma$, there exists a holomorphic function $\psi$ such that
\begin{equation} 
\label{eqn:secGauss}
g_{m} = (1+\psi\bar{\psi})e_{0} + (\psi +\bar{\psi})e_{1} - i(\psi - \bar{\psi}) e_{2} + (-1+\psi\bar{\psi})e_{3}
\end{equation}
is a lift of $G_{m}$. We then define a holomorphic 1-form $\eta$ so that $q = d\psi \eta + d\bar{\psi}\bar{\eta}$. The closed 1-form $\zeta_{m}= Ad_{T_{m}}\zeta$ is then identified with 
\[ \xi_{m}:= \begin{pmatrix} -\psi & \psi^{2} \\ -1 &\psi\end{pmatrix}\eta.\]
Now parallel sections $x$ of $d+m\zeta$ are of the form $x=T^{-1}_{m}\mathfrak{c} = T^{m}_{-m}\mathfrak{c}$, using~\eqref{eqn:iterate}. $T^{m}_{-m}$ satisfies $dT^{m}_{-m}= -mT^{m}_{-m}\zeta_{m}$ and correspondingly $\Psi_{m}= F^{-1}_{m}$ satisfies 
\begin{equation}
\label{eqn:secondrep}
d\Psi_{m} = -m\Psi_{m}\xi_{m} = m\Psi_{m}\begin{pmatrix} \psi & -\psi^{2} \\ 1 &-\psi\end{pmatrix}\eta.
\end{equation}
This formulation shows that the representation~\eqref{eqn:hyprep} coincides with the representation of~\cite{UY1992} when $\mu=-1$ and with the representation of~\cite{AA1998} when $\mu=1$. The duality between~\eqref{eqn:Psi} and~\eqref{eqn:secondrep} has been remarked upon in~\cite{BPS2003,RUY1997}. 

\subsection{Linear Weingarten surfaces of Bryant type}
A surface $x:\Sigma\to \mathbb{H}^{3}$ is a \textit{linear Weingarten surface of Bryant type} if the mean curvature $H$ and Gauss curvature $K$ of $x$ satisfy a relation
\begin{equation} 
\label{eqn:brytype}
(\mu+1)K-2\mu H + \mu-1=0
\end{equation}
for some $\mu\in\mathbb{R}$. A Weierstrass-type representation was developed for these surfaces in~\cite{GMM2004}. The middle sphere congruence of such surfaces is given by $x^{M} = x-\frac{\mu+1}{2}\tilde{g}$, where $\tilde{g}\in \Gamma G$ such that $(\tilde{g},x)=-1$. Three cases emerge\footnote{This analysis is analogous to that performed in~\cite[Section 4.6]{BHR2014} for parallel families of linear Weingarten surfaces in hyperbolic space.}:
\begin{itemize}
\item if $\mu<0$ then $x^{M}$ is a CMC-$\frac{1}{\sqrt{-\mu}}$ surface in $\mathbb{H}^{3}(\frac{1}{\mu})$, 
\item if $\mu>0$ then $x^{M}$ is CMC-$\frac{1}{\sqrt{\mu}}$ surface in $\mathbb{S}^{2,1}(\frac{1}{\mu})$, 
\item if $\mu=0$ then $x^{M}$ is an intrinsically flat surface in $\mathcal{L}$.
\end{itemize}
We thus have that 
\[x^{M} = \Psi_{m}\begin{pmatrix} 1 & 0\\ 0&-\mu\end{pmatrix}\Psi^{*}_{m}\]
where $\Psi_{m}:\Sigma\to SL(2,\mathbb{C})$ satisfies~\eqref{eqn:secondrep}. 
Using the Hermitian model, the lift $g_{m}$ from~\eqref{eqn:secGauss} of the secondary Gauss map is given by
\[ g_{m} = 2\begin{pmatrix}
|\psi|^{2} & \psi\\ \bar{\psi} & 1
\end{pmatrix}.\]
Now, since $\tilde{g}\in \Gamma G$ satisfies $(\tilde{g},x)=-1$ and $T^{-1}_{m}g_{m}\in \Gamma G$, 
\[
\tilde{g} =- \frac{T_{m}^{-1}g_{m}}{(T_{m}^{-1}g_{m},x)} =- \frac{ \Psi_{m}g_{m}\Psi_{m}^{*}}{\left(\Psi_{m}g_{m}\Psi_{m}^{*}, \Psi_{m}\begin{pmatrix} 1 & 0\\ 0&-\mu\end{pmatrix}\Psi^{*}_{m}\right)}= \frac{\Psi_{m}g_{m}\Psi_{m}^{*}}{1-\mu|\psi|^{2}}.
\]
Thus\footnote{Note that in order for this expression to be well defined, one must assume that $1-\mu|\psi|^{2}$ is nowhere zero. This coincides with the assumption of~\cite{GMM2004}.}, 
\begin{align*}
 x = x^{M} + \frac{\mu+1}{2}\tilde{g} &= \Psi_{m}\left( \begin{pmatrix} 1 & 0\\ 0&-\mu\end{pmatrix} + \frac{\mu+1}{1-\mu|\psi|^{2}} \begin{pmatrix}
|\psi|^{2} & \psi \\ \bar{\psi} & 1
\end{pmatrix} \right) \Psi^{*}_{m}\\
&= \frac{1}{1-\mu|\psi|^{2}} \Psi_{m}\begin{pmatrix} 1+|\psi|^{2} & (\mu+1)\psi \\ (\mu+1)\bar{\psi} & 1+\mu^{2}|\psi|^{2}\end{pmatrix}\Psi_{m}^{*}. 
\end{align*}
By setting $H:= \Psi_{m}\begin{pmatrix} i\psi & i\\i & 0\end{pmatrix}$ we have that
\[ x = H\begin{pmatrix} \frac{1+\mu^{2}|\psi|^{2}}{1-\mu|\psi|^{2}} & \mu\bar{\psi}\\ \mu \psi & 1-\mu |\psi|^{2}\end{pmatrix}H^{*}.\]
We deduce from~\eqref{eqn:secondrep} that 
\[ H^{-1}dH= \begin{pmatrix} 0 & m \eta\\ d\psi & 0\end{pmatrix}.\]
Hence, we have recovered the representation of~\cite{GMM2004}. Moreover, in the case that $\mu=0$, we obtain the representation of~\cite{GMM2000} for flat surfaces in $\mathbb{H}^{3}$ . 

\begin{remark}
\textit{Linear Weingarten surfaces of Bianchi type}, that is surfaces satisfying~\eqref{eqn:brytype} in $\mathbb{S}^{2,1}$, were shown to admit a Weierstrass-type representation in~\cite{AE2007}. An analogous analysis can be performed for these surfaces as above. 
\end{remark}

\bibliographystyle{plain}
\bibliography{wtype}

\end{document}